\documentclass[10pt,reqno]{amsart}

\usepackage[text={160mm,240mm},centering]{geometry}            
\usepackage{amsmath}
\usepackage{comment}
\usepackage{graphicx}
\usepackage{amssymb}
\usepackage{dsfont}
\usepackage{tikz-cd}
\usepackage{amsmath,amsthm}
\usepackage{mathabx}
\usepackage[mathscr]{eucal}
\usepackage[all]{xy}
\usepackage{mathrsfs}
\usepackage{hyperref}
\usepackage{color}
\definecolor{ccr}{RGB}{10,110,150}
\hypersetup{colorlinks=true,linkcolor=ccr,anchorcolor=ccr,citecolor=ccr,urlcolor=ccr}

\newcommand{\m}{\mathrm}
\newtheorem{theorem}{Theorem}[section]
\newtheorem{lemma}[theorem]{Lemma}
\numberwithin{equation}{section}

\newtheorem{conjecture}[theorem]{Conjecture}
\newtheorem{corollary}[theorem]{Corollary}
\newtheorem{proposition}[theorem]{Proposition}

\newtheorem{question}[theorem]{Question}
\newtheorem{definition-lemma}[theorem]{Definition}
\newtheorem{definition-theorem}[theorem]{Definition}

\usepackage{enumitem}

\theoremstyle{definition}
\newtheorem{example}[theorem]{Example}
\newtheorem{definition}[theorem]{Definition}

\usepackage{bm}

\newtheorem{remark}[theorem]{Remark}

\newtheorem*{ack}{Acknowledgements}

\newtheorem*{ai}{AI Declaration}
\newtheorem*{org}{Organization of the paper}

\newcommand{\hfk}{\widehat{\mathrm{HFK}}}
\newcommand{\ccfl}{\mathcal{CFK}^-}
\newcommand{\cfk}{\mathrm{CFK}}
\newcommand{\ch}{V}
\newcommand{\hfp}{\mathrm{HFK}'}
\allowdisplaybreaks

\title{Torsion order and the next-to-top term of knot Floer homology}

\author{Kai Zhu}
\address{Beijing International Center for Mathematical Research, Peking University, Beijing 100871,  P. R. China.}
\email{kzhu25@stu.pku.edu.cn}

\begin{document}

	\renewcommand{\abstractname}{Abstract}
	\begin{abstract}
    Sivek conjectured that, for every nontrivial knot $K\subset S^3$, the dimension of the next-to-top term of $\hfk(K)$ is at least that of the top term. We show that this inequality holds for knots with $U$-torsion order one in either the minus theory $\m{HFK}^-$ or the unoriented theory $\m{HFK}'$. Combined with Eftekhary's work, this implies that knots with rational unknotting number one satisfy the conjecture. In particular, the result applies to every nontrivial satellite knot whose pattern can be unknotted in the solid torus via a single rational tangle replacement, such as the $(2,2n+1)$-cable of a nontrivial knot for any integer $n$. In the appendix, we study the possible pairs of the two torsion orders of a knot.
	\end{abstract}

	\maketitle
	
	\parskip=.3em
	
	\setcounter{page}{1}

\section{Introduction}

Knot Floer homology was introduced independently by Ozsv\'ath and Szab\'o~\cite{OzsvathSzabo2004} and by Rasmussen~\cite{Rasmussen2003} for null-homologous knots in closed, oriented $3$-manifolds. For a knot $K\subset S^3$, the hat version $\hfk(K)$ is a finite-dimensional vector space over $\mathbb F=\mathbb F_2$, with Maslov grading $m$ and Alexander grading $a$. Thus
\[
    \hfk(K)=\bigoplus_{m,a\in\mathbb{Z}}\hfk_m(K,a).
\]

Let $g(K)$ denote the Seifert genus of $K$. The top term $\hfk(K,g(K))$ is always nonzero~\cite{OzsvathSzabo2004b}, and it has dimension one precisely when $K$ is fibered~\cite{Ni2007}. Hence it is natural to ask what information the other terms of $\hfk(K)$ contain. But before we deduce some interesting topological applications, we need to know whether these terms are nontrivial.

Baldwin and Vela-Vick~\cite[Question~1.12]{BaldwinVelaVick2018} asked whether the next-to-top term $\hfk(K,g(K)-1)$ is nontrivial for every knot of positive genus. More specifically, Sivek~\cite[Question~1.13]{BaldwinVelaVick2018} proposed the following conjecture.
	
\begin{conjecture}\label{conjsivek}
    For every nontrivial knot $K\subset S^3$,
    \[
        \dim\hfk(K,g(K)-1)\ge \dim\hfk(K,g(K)).
    \]
\end{conjecture}

Known cases of the conjectured inequality and its extensions to other $3$-manifolds include
\begin{itemize}
    \item nontrivial knots with thin knot Floer homology, by the decomposition in~\cite[Lemma~7]{Petkova2013},
    \item nontrivial $L$-space knots~\cite{HeddenWatson2018},
    \item fibered knots of positive genus in closed, oriented $3$-manifolds~\cite[Theorem~1.1]{BaldwinVelaVick2018},
    \item  with rational coefficients, null-homologous knots in generalized \(L\)-spaces with \(b_1\le 1\) whose top term is supported in a single \(\mathbb{Z}/2\mathbb{Z}\)-Maslov grading ~\cite[Theorem~1.1]{Ni2022},
    \item nontrivial knots that are closures of $3$-braids~\cite[Corollary~1.3]{Chen2025},
    \item nontrivial knots with $w$-avoiding exterior~\cite[Theorem~1.3]{BinnsWan2026}.
\end{itemize}

In this note, we establish Conjecture~\ref{conjsivek} for several new families of knots using the oriented and unoriented torsion orders in knot Floer homology. Recall that for a knot $K\subset S^3$, the minus version $\m{HFK}^-(K)$ and the unoriented version $\m{HFK}'(K)$ are finitely generated $\mathbb{F}[U]$-modules. They are obtained by taking the homology of the full knot Floer complex over $\mathbb{F}[U,V]$ after setting $V=0$ and $V=U$, respectively.

The oriented torsion order of $K$~\cite{JuhaszMillerZemke2020} is defined by
\[
    \m{Ord}_U(K)=\min\{n\in\mathbb{Z}_{\ge 0}: U^n\cdot \m{Tor}(\m{HFK}^-(K))=0\},
\]
and the unoriented torsion order~\cite{GM2023} by
\[
    \m{Ord}'_U(K)=\min\{n\in\mathbb{Z}_{\ge 0}: U^n\cdot \m{Tor}(\m{HFK}'(K))=0\}.
\]

\begin{theorem}\label{ord}
    Let $K\subset S^3$ be a knot with $\m{Ord}_U(K)=1$. Then for every integer $m$,
    \begin{equation}\label{2inequal}
        \dim\hfk_{m-1}(K,g(K)-1)\ge 2\dim\hfk_m(K,g(K))+c_{m,K},
    \end{equation}
    where $c_{m,K}\in\{-1,0,1\}$ is given by~\eqref{correctionterm}. The third-to-top term satisfies
    \[
        \dim\hfk_{m-2}(K,g(K)-2)\ge\dim\hfk_m(K,g(K)).
    \]
\end{theorem}

\begin{remark}
For a related nonvanishing result of the third-to-top term of knot Floer homology, Truong~\cite[Theorem~1.5]{Truong2026} showed that the term is nonzero for a fibered knot whose monodromy is neither left-veering nor right-veering.
\end{remark}

\begin{corollary}\label{unknotting1}
    Let $K\subset S^3$ be a knot with unknotting number $\m{u}(K)=1$. Then for every integer $m$,
    \[
        \dim \hfk_{m-1}(K,g(K)-1)\ge \dim \hfk_m(K,g(K)).
    \]
\end{corollary}
\begin{proof}
    By~\cite[Proposition~4.1]{AE20}, $\m{Ord}_U(K)\le\m{u}(K)=1$. If the order were zero, $\m{HFK}^-(K)$ would be a single free tower. Then~\eqref{leq} would give $\dim\hfk(K)=1$, forcing $K$ to be the unknot~\cite{OzsvathSzabo2004b,OzsvathSzabo2004}. Thus $\m{Ord}_U(K)=1$. Now Theorem~\ref{ord} gives a lower bound of $\max\{2\dim\hfk_m(K,g(K))-1,0\}\ge \dim\hfk_m(K,g(K))$ for the next-to-top dimension in Maslov grading $m-1$.
\end{proof}

\begin{remark}
    Bhat, Li, and Ye~\cite[Corollary~1.11]{BhatLiYe2025} proved that the next-to-top term of instanton knot homology $\m{KHI}(S^3,K)$ is nonzero when $\m{u}(K)=1$. They suggested that their argument might extend to knot Floer homology via the unoriented theory~\cite[Remark~1.12]{BhatLiYe2025}. Corollary~\ref{unknotting1} gives a direct proof of the corresponding nonvanishing result using the minus theory.
\end{remark}

\begin{theorem}\label{tor}
    Let $K\subset S^3$ be a knot with $\m{Ord}'_U(K)=1$. Then for every integer $m$,
    \[
        \dim \hfk_{m-1}(K,g(K)-1)\ge \dim \hfk_m(K,g(K)).
    \]
\end{theorem}

Eftekhary's work~\cite{eftekhary2026tangle} gives a family of knots with unoriented torsion order one. We recall the relevant definitions from~\cite{iltgen2025khovanov}.

\begin{definition}
    Two knots $K$ and $K'$ are related by a \textit{rational tangle replacement} if $K'$ can
be obtained from $K$ by replacing a rational tangle $T$ in $K$ with another rational
tangle $T'$ in the same ball with the same four endpoints. If the arcs of $T$ and $T'$ induce the same pairing of the endpoints, the
replacement is \textit{proper}.
\end{definition}

\begin{definition}
    The \textit{rational unknotting number} $\m{u}_r(K)$ is the minimum number of rational tangle replacements needed to turn $K$ into the unknot.
The \textit{proper rational unknotting number} $\m{u}_{pr}(K)$ is the corresponding minimum when all replacements are proper.
\end{definition}

\begin{figure}[htbp]
    \centering
    \includegraphics[width=0.54\linewidth]{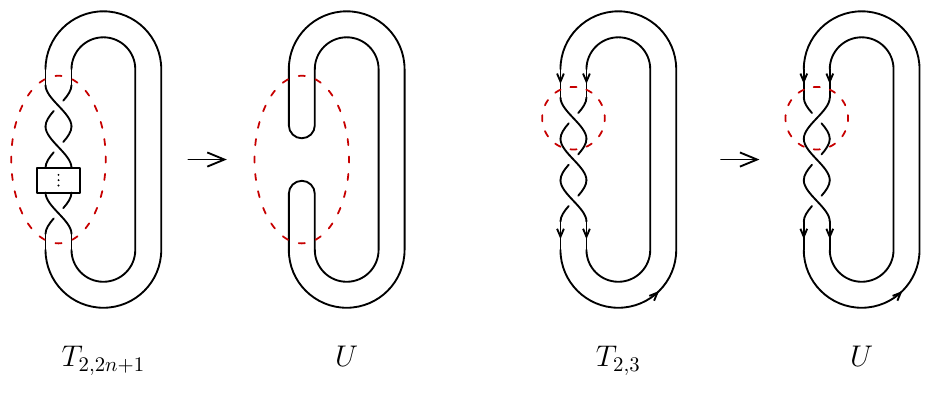}
    \caption{The rational tangle replacement on the left gives $\m{u}_{r}(T_{2,2n+1})=1$ for $|2n+1|>1$. The proper replacement on the right gives $\m{u}_{pr}(T_{2,3})=1$.}
    \label{fig:tangle-replacement}
\end{figure}

\begin{remark}\label{uuu}
    A crossing change is a proper rational tangle replacement. Thus $\m{u}_r(K)\le \m{u}_{pr}(K)\le \m{u}(K)$, where $\m{u}(K)$ is the unknotting number. Figure~\ref{fig:tangle-replacement} illustrates both types of replacement.
\end{remark}

\begin{proposition}\label{3}
    Let $K\subset S^3$ be a knot with rational unknotting number $\m{u}_r(K)=1$. Then for every integer $m$,
    \[
        \dim \hfk_{m-1}(K,g(K)-1)\ge \dim \hfk_m(K,g(K)).
    \]
\end{proposition}

\begin{proof}
In Eftekhary's terminology~\cite{eftekhary2026tangle}, an $n$-tangle replacement takes place in a ball $B$ for which $T=B\cap K$ is a trivial $n$-strand tangle. It replaces $T$ by $f(T)$ for a diffeomorphism $f:B\to B$ that preserves the set of endpoints. Any two rational tangles with the same endpoints are related by such a diffeomorphism, so a rational tangle replacement is a $2$-tangle replacement.

Eftekhary defines the tangle distance $\m{d}_t(K,K')$ as the minimum $n$ such that an $(n+1)$-tangle replacement changes $K$ into $K'$. By \cite[Theorem~1.1]{eftekhary2026tangle},
\begin{equation}
    \m{d}_t(K,K')\ge \m{Ord}'_U(K)-\m{Ord}'_U(K').
\end{equation}
If $\m{u}_r(K)=1$, a single $2$-tangle replacement turns $K$ into the unknot $\m U$. Since $\m{Ord}'_U(\m U)=0$, we obtain $\m{Ord}'_U(K)\le 1$. If the order were zero, $\m{HFK}'(K)$ would consist of a single free tower, so reduction modulo $U$ would give $\dim\hfk(K)=1$ by the exact sequence~\eqref{pleq}. This implies  $K$ is unknot from ~\cite{OzsvathSzabo2004b} which is a contradiction. Thus $\m{Ord}'_U(K)=1$, and the result follows from Theorem~\ref{tor}.
\end{proof}

We now give examples of knots with rational unknotting number one. We call a pattern $P\subset S^1\times\mathbb D^2$ \emph{ rational unknotted} if a single rational tangle replacement in the solid torus turns $P$ into a knot that bounds an embedded disk there.

\begin{remark}
    Dai, Hedden, Mallick, and Stoffregen call this a pattern of rational unknotting number one~\cite[Definition~2.2]{daieftal2024rank}. However, one should be careful that this is not the same as $\m{u}_r(P)=1$: the latter means that it can be isotoped to an unknot in $S^3$, not merely in the solid torus. Hence we use different notation to avoid confusion.
\end{remark}

For a  rational unknotted pattern $P$, the same replacement in $P(K)$ gives a knot that bounds a disk in a tubular neighborhood of $K$. Hence $\m{u}_r(P(K))\le 1$, with equality whenever $P(K)$ is nontrivial. Proposition~\ref{3} gives the following corollary.

\begin{corollary}\label{5}
    Let $K\subset S^3$ be a knot and $P$ a  rational unknotted pattern. If $P(K)$ is nontrivial, then for every integer $m$,
    \[
        \dim \hfk_{m-1}\bigl(P(K),g(P(K))-1\bigr)\ge \dim \hfk_m\bigl(P(K),g(P(K))\bigr).
    \]
\end{corollary}

\begin{figure}[htbp]
    \centering
    \includegraphics[width=0.99\linewidth]{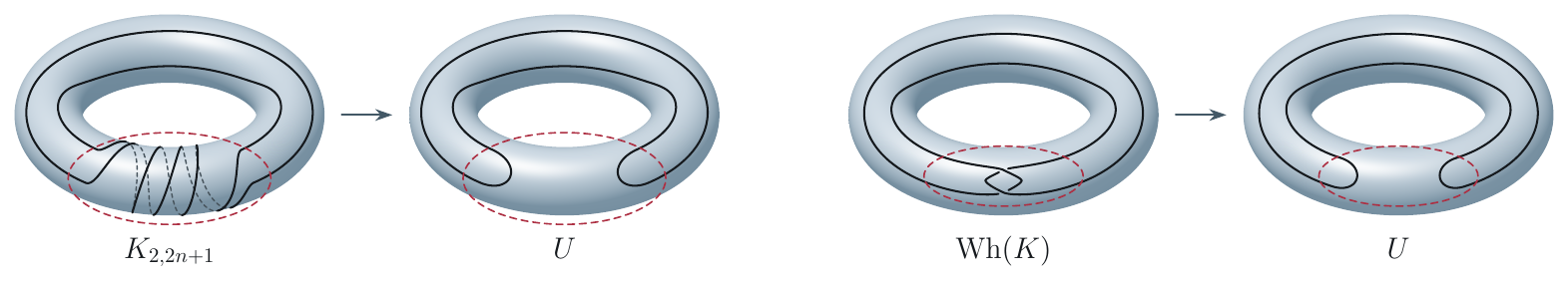}
    \caption{Rational tangle replacements unknotting the $(2,2n+1)$-cable pattern (left) and the Whitehead double pattern (right) in the solid torus.}
    \label{fig:solid-torus-patterns}
\end{figure}

\begin{example}
    Rational tangle patterns give simple examples. Such a pattern is the horizontal closure of a rational tangle in the solid torus, with the closure taken to be a knot. Replacing the tangle by the standard $\infty$-tangle gives a knot bounding a disk in the solid torus. Figure~\ref{fig:solid-torus-patterns} illustrates the replacement for a $(2,2n+1)$-cable pattern and a Whitehead double pattern. Further examples arise from Conway tangles with unknotted horizontal closure and from strongly invertible knots~\cite[Section~2.4]{daieftal2024rank}.
\end{example}

\begin{remark}
From the preceding example, Corollary~\ref{5} applies to nontrivial Whitehead doubles and $(2,2n+1)$-cables. For Whitehead doubles, the inequality also follows from Hedden's computation~\cite{hedden2007knot}. For a fixed knot $K$ and all sufficiently large positive integers $n$. Hedden's cabling formula~\cite[Corollary~1.3]{Hedden2005} gives
 \[\hfk_{m-1}(K_{2,2n+1},g(K_{2,2n+1})-1)\cong\hfk_m(K_{2,2n+1},g(K_{2,2n+1}))\cong\hfk_m(K,g(K))\]
    for every $m$. Thus the top two terms have equal dimensions in this range. Corollary~\ref{5} establishes the dimension inequality for every integer $n$, provided that $K_{2,2n+1}$ is nontrivial.  
\end{remark}

\begin{org}
Section~\ref{sec:prelim} recalls the basic background on knot Floer homology. In Section~\ref{sec:oriented}, we study the action of \(U\) on the top Alexander grading of \(\m{HFK}^-\) and use the staircase and box decompositions of the knot Floer chain complex for knots of torsion order one to prove Theorem~\ref{ord}. Section~\ref{sec:unoriented} proves Theorem~\ref{tor} and in fact provides a unified proof that knots with either oriented or unoriented torsion order one satisfy Conjecture~\ref{conjsivek}. Appendix~\ref{sec:torsion-orders} compares the two torsion orders and explores which pairs of integers can be realized as the oriented and unoriented torsion orders of a knot.
\end{org}

\begin{ack}
The author thanks his advisor, Fan Ye, for proposing this project, suggesting the use of unoriented knot Floer homology, sharing ChatGPT conversation records concerning the unoriented torsion order, and pointing out Proposition~\ref{stru} as well as the idea for a simplified proof of Proposition~\ref{prop4.3}.
\end{ack}

\begin{ai}
The author independently developed the original argument for the oriented torsion order, including Theorem~\ref{ord}, Corollary~\ref{unknotting1}, and Theorem~\ref{waction}. ChatGPT 6 Astra supplied the original proof of Theorem~\ref{tor} in the unoriented setting, as presented in the first version of this paper; the proof appearing in this version has been simplified by the author. ChatGPT 6 Astra suggested the example of \((2,2n+1)\)-cables in Corollary~\ref{5}, which motivated the author to consider the more general class of rational unknotted patterns. Appendix~\ref{sec:torsion-orders} grew out of discussions between the author and ChatGPT. The author further used ChatGPT 6 Astra to create figures and edit the English. The author takes full responsibility for the mathematical arguments and the final text.    
\end{ai}

\section{Preliminaries on knot Floer homology}\label{sec:prelim}

In this section, we recall the constructions in knot Floer homology that define the oriented and unoriented torsion orders.

Let $(\Sigma,\bm{\alpha},\bm{\beta},w,z)$ be an admissible doubly pointed Heegaard diagram for a knot $K\subset S^3$. The full knot Floer chain complex $\mathcal{CFK}^-(K)$ is the free $\mathbb{F}[U,V]$-module generated by intersection points between the associated tori $\mathbb{T}_{\bm{\alpha}}$ and $\mathbb{T}_{\bm{\beta}}$ in $\m{Sym}^h(\Sigma)$, where $h$ is the genus of $\Sigma$. The differential $\partial :\mathcal{CFK}^-(K)\to \mathcal{CFK}^-(K)$ is the $\mathbb{F}[U,V]$-linear map defined on generators by \begin{equation}
    \partial\mathbf{x}=\sum_{\mathbf{y}\in \mathbb{T}_{\bm{\alpha}}\cap \mathbb{T}_{\bm{\beta}}} \sum_{\substack{\phi\in \pi_2(\mathbf{x},\mathbf{y}) \\ \mu(\phi)=1}} \#\widehat{\mathcal{M}}(\phi)\cdot U^{n_w(\phi)}V^{n_z(\phi)}\cdot \mathbf{y}.
\end{equation}

Here $\#\widehat{\mathcal M}(\phi)$ is the mod-two count of unparametrized holomorphic disks in the class $\phi$. The complex $\ccfl(K)$ has a bigrading $(\m{gr}_w,\m{gr}_z)$ and the Alexander grading $A=\dfrac{1}{2}\bigl( \m{gr}_w-\m{gr}_z\bigr)$. Multiplication by $U$ and $V$ has grading shifts $(-2,0)$ and $(0,-2)$, respectively, while $\partial$ has grading shift $(-1,-1)$.

The full knot Floer complex has various algebraic specializations. Setting $V=0$, we obtain the minus complex \cite{OzsvathSzabo2004},
\[
    \cfk^-(K)=\ccfl(K)/(V),\quad \partial^-\mathbf{x}=\sum_{\mathbf{y}\in \mathbb{T}_{\bm{\alpha}}\cap \mathbb{T}_{\bm{\beta}}} \sum_{\substack{\phi\in \pi_2(\mathbf{x},\mathbf{y}) \\ \mu(\phi)=1,n_z(\phi)=0}} \#\widehat{\mathcal{M}}(\phi)\cdot U^{n_w(\phi)} \mathbf{y}.
\]
This is a free $\mathbb{F}[U]$-module. Its homology is $\m{HFK}^-(K)$ and is bigraded by $(M,A)$, where $M=\m{gr}_w$ is the Maslov grading and $A$ is the Alexander grading.

Setting $U=V$, we obtain the unoriented knot Floer complex \cite{OzsvathStipsiczSzabo2017},
\[
    \cfk'(K)=\ccfl(K)/(U-V).
\]
We regard it as a module over $\mathbb{F}[U]$, using $U$ for the common image of the two variables. With this convention the differential is
\[
    \partial'\mathbf{x}=\sum_{\mathbf{y}\in \mathbb{T}_{\bm{\alpha}}\cap \mathbb{T}_{\bm{\beta}}} \sum_{\substack{\phi\in \pi_2(\mathbf{x},\mathbf{y}) \\ \mu(\phi)=1}} \#\widehat{\mathcal{M}}(\phi)\cdot U^{n_w(\phi)+n_z(\phi)}\cdot \mathbf{y}.
\]
Its homology is $\m{HFK}'(K)$. It is graded by $\delta=\dfrac{1}{2}(\m{gr}_w+\m{gr}_z)=M-A$, and both $\partial'$ and multiplication by $U$ lower $\delta$ by $1$. The Alexander grading does not descend through the relation $U=V$.

Finally, setting $U=V=0$, we obtain the hat knot Floer complex $\widehat{\cfk}(K)$ with homology $\hfk(K)$. Equivalently, setting $U=0$ in either $\cfk^-(K)$ or $\cfk'(K)$ gives $\widehat{\cfk}(K)$. Hence $\hfk(K)$ admits both bigradings $(M,A)$ and $(\delta,A)$, related by $\delta=M-A$.

Since $K$ is a knot in $S^3$, $\m{HFK}^-(K)$ and $\m{HFK}'(K)$ each contain a single free tower \cite{OzsvathSzabo2004, OzsvathStipsiczSzabo2017}. By the graded structure theorem for finitely generated modules over the principal ideal domain $\mathbb{F}[U]$, we have \begin{align}
    \m{HFK}^-(K)&\cong \mathbb{F}[U]_{(m,a)}\oplus\bigoplus_{i} \bigl(\mathbb{F}[U]/(U^{n_i})\bigr)_{(m_i,a_i)}, \\
    \m{HFK}'(K)&\cong \mathbb{F}[U]_{\delta_0}\oplus\bigoplus_{j} \bigl(\mathbb{F}[U]/(U^{\ell_j})\bigr)_{\delta_j}.
\end{align}
The subscripts give the $(M,A)$-bigradings of the generators in the first decomposition and their $\delta$-gradings in the second. All torsion exponents are positive integers. It follows that the torsion order of $K$ is
\[
    \m{Ord}_U(K)=\min\{n\in\mathbb{Z}_{\ge 0}:U^n\cdot \m{Tor}\bigl( \m{HFK}^-(K)\bigr)=0 \}=\max\bigl(\{0\}\cup\{n_i\}_i\bigr).
\]
The unoriented torsion order is
\[
    \m{Ord}_U'(K)=\min\{n\in\mathbb{Z}_{\ge 0}:U^n\cdot \m{Tor}\bigl( \m{HFK}'(K)\bigr)=0 \}=\max\bigl(\{0\}\cup\{\ell_j\}_j\bigr).
\]

\section{Knots with oriented torsion order one}\label{sec:oriented}

We first prove Conjecture~\ref{conjsivek} for knots of oriented torsion order one by examining the action of $U$ in the top Alexander grading. We then use the structure of their knot Floer complexes to prove  a stronger result, namely Theorem \ref{ord}.

\subsection{The action of $U$ on the top Alexander grading}

\begin{theorem}\label{waction}
    Let $K\subset S^3$ be a knot with Seifert genus $g$, and fix an integer $m$. If $U\cdot \m{HFK}^-_m(K,g)=0$, then
\[
    \dim \hfk_{m-1}(K,g-1)\ge \dim \hfk_m(K,g).
\]
\end{theorem}

\begin{proof}
The short exact sequence
\[
    0\to \cfk^-(K)\overset{U}{\to}\cfk^-(K)\to \widehat{\cfk}(K)\to 0
\]
induces the long exact sequence
\begin{equation}\label{leq}
    \begin{aligned}
        \cdots\to \m{HFK}^-_{m+1}(K,a)\overset{U}{\to} \m{HFK}^-_{m-1}(K,a-1)&\to \hfk_{m-1}(K,a-1)\\
        &\overset{\gamma}{\to} \m{HFK}^-_{m}(K,a)\overset{U}{\to} \m{HFK}^-_{m-2}(K,a-1)\to \cdots.
    \end{aligned}
\end{equation}
Indeed, if $x$ lifts a cycle $\bar x$ and $\partial^-x=Uy$, then $\gamma([\bar x])=[y]$. Since $\partial^-$ and $U$ have bidegrees $(-1,0)$ and $(-2,-1)$, the connecting homomorphism $\gamma$ has bidegree $(1,1)$.

At $a=g+1$, the minus groups in Alexander grading $g+1$ vanish. Thus~\eqref{leq} identifies $\m{HFK}^-_m(K,g)$ with $\hfk_m(K,g)$. At $a=g$, the hypothesis on $U$ makes $\gamma$ surjective, which gives the inequality.
\end{proof}

\begin{corollary}\label{tor1}
    Let $K\subset S^3$ be a knot of Seifert genus $g$ with $\m{Ord}_U(K)=1$. Then for every integer $m$,
\[
    \dim \hfk_{m-1}(K,g-1)\ge \dim \hfk_m(K,g).
\]
\end{corollary}

\begin{proof}
Recall $-\tau(K)$ is  the Alexander grading of the generator of the free tower of $\m{HFK}^-(K)$ \cite[Proposition 6.1.4]{OzsvathStipsiczSzabo2015}. If $-\tau(K)<g$, the entire top group is torsion and is killed by $U$, so Theorem~\ref{waction} applies.

If $-\tau(K)=g$, consider the mirror $\m{m}(K)$. Its minus complex is the $\mathbb F[U]$-dual of that of $K$, so the torsion exponents are unchanged and its free tower starts in Alexander grading $-g$. Since $K$ is nontrivial, $g>0$, and the preceding case applies to $\m{m} (K)$. Mirror duality and conjugation symmetry~\cite{OzsvathSzabo2004} give
\[
    \dim\hfk_j(\m{m}( K),a)=\dim\hfk_{2a-j}(K,a).
\]
Taking $j=2g-m$ in the inequality for $\overline K$ proves the claim.
\end{proof}

\begin{remark}
The hypothesis of Theorem~\ref{waction} also holds for knots of arbitrarily large torsion order. For $n\ge2$, the first step in the staircase of the positive torus knot $T_{n,n+1}$ has length one, so $U$ kills its top group. On the other hand, $\m{Ord}_U(T_{n,n+1})=n-1$~\cite[Corollary~5.3]{JuhaszMillerZemke2020}. Thus, while the proof of Theorem~\ref{waction} is relatively simple, it may yield further families of knots that satisfy Conjecture~\ref{conjsivek}.
\end{remark}

\subsection{Proof of Theorem~\ref{ord}}\label{subsec:boxes}

Set $\mathcal R=\mathbb F[U,V]/(UV)$ and
\[
    \mathcal{CFK}_{\mathcal R}(K)=\mathcal{CFK}^-(K)\otimes_{\mathbb F[U,V]}\mathcal R.
\]
Its reduction modulo $(U,V)$ is the hat knot Floer complex. We use the following description of $\mathcal{CFK}_{\mathcal R}(K)$.

\begin{proposition}[{\cite[Lemma~2.8]{CZZ2026}, cf.~\cite[Lemma~7]{Petkova2013}}]\label{standardord1}
    If $\m{Ord}_U(K)=1$, then there is a bigraded chain homotopy equivalence
    \[
        \mathcal{CFK}_{\mathcal R}(K)\simeq\mathcal C_{\tau(K)}\oplus\mathcal B_1\oplus\cdots\oplus\mathcal B_k,
    \]
    where $k\ge0$, $\mathcal C_{\tau(K)}$ is a staircase with unit steps, and each $\mathcal B_i$ is a length-one box. The staircase $\mathcal C_0$ has a single generator.
\end{proposition}

\begin{figure}[htbp]
    \centering
    \includegraphics[width=0.7\linewidth]{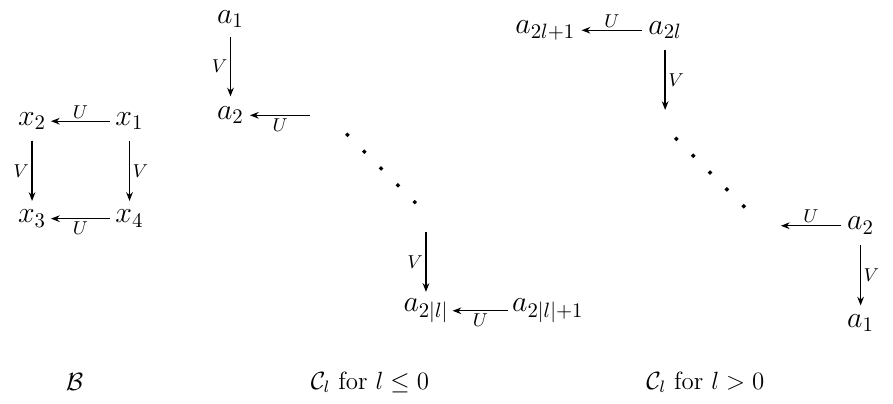}
    \caption{A length-one box and the two staircase types for $\mathcal{CFK}_{\mathcal R}(K)$. An arrow from $x$ to $y$ labeled $U$ or $V$ indicates the corresponding term in $\partial x$. Compare~\cite[Figure~2]{Petkova2013}.}
    \label{fig:oriented-models}
\end{figure}

\begin{proof}[Proof of Theorem~\ref{ord}]
 Every differential in the model of Proposition~\ref{standardord1} vanishes modulo $(U,V)$, so its generators give a bigraded basis of $\hfk(K)$.

Consider a box that meets Alexander grading $g(K)$. Only its top generator $x_2$ can lie in this grading. If $(M,A)(x_2)=(m,g)$, the unit arrows in Figure~\ref{fig:oriented-models} place $x_1,x_3$ in $(m-1,g-1)$ and $x_4$ in $(m-2,g-2)$. Thus every box generator in $(m,g)$ contributes two generators one level below and one generator two levels below, all on the same $M-A$ diagonal.

The staircase has one generator in each bidegree \cite[Proposition 6.1.4]{OzsvathStipsiczSzabo2015},
\[
    (M,A)=(a-\tau(K),a),\qquad -|\tau(K)|\le a\le |\tau(K)|.
\]
Indeed, its unit steps preserve $M-A$, and the endpoint representing the generator of $\widehat{\m{HF}}(S^3)$ after setting $U=0$, $V=1$ has $(M,A)=(0,\tau(K))$. The other endpoint has Alexander grading $-\tau(K)$ by basepoint symmetry. If $|\tau(K)|=g(K)$, the staircase contributes once at both $g(K)$ and $g(K)-1$. If $|\tau(K)|=g(K)-1$, it contributes only at $g(K)-1$ among these two gradings. In either case, the contribution lies on the diagonal $m-g(K)=-\tau(K)$. For $|\tau(K)|\le g(K)-2$ it contributes at neither grading. Adding the box contributions gives~\eqref{2inequal} with
\begin{equation}\label{correctionterm}
    c_{m,K}=\begin{cases}
       -1 & \text{if }m=g(K)-\tau(K)\text{ and }|\tau(K)|=g(K),\\
        1 & \text{if }m=g(K)-\tau(K)\text{ and }|\tau(K)|=g(K)-1,\\
        0 & \text{otherwise}.
    \end{cases}
\end{equation}

Finally, if the staircase reaches $g(K)$, then $|\tau(K)|=g(K)\ge1$, so it also reaches $g(K)-2$. Together with the bottom generators of the boxes above, this gives
\[
    \dim\hfk_{m-2}(K,g(K)-2)\ge\dim\hfk_m(K,g(K)),
\]
which implies the second inequality of Theorem~\ref{ord}.
\end{proof}

\section{Knots with unoriented torsion order one}\label{sec:unoriented}

In this section, we give a unified proof of Theorem~\ref{tor} and Corollary~\ref{tor1} in a more general algebraic framework.

Let $(C,\partial)$ be a finitely generated free complex over $\mathbb F[U]$, graded by $\delta$, with both $\partial$ and multiplication by $U$ of degree $-1$. Write $\widehat{C}=C/UC$ and $\widehat{V}=\m{H}\bigl(\widehat{C},\widehat\partial\bigr)$. For the Alexander grading estimates below, assume that $\widehat{V}$ has a further grading $\widehat{V}=\bigoplus_{\delta,a}\widehat{V}_{\delta,a}$ with integer Alexander degrees $-g\le a\le g$, where $g$ is a positive integer. Assume that $\dim \widehat{V}_{\delta,a}=\dim \widehat{V}_{\delta,-a}$ for all $\delta,a$.

The following normal form is the one-variable version of the simplified bases used in knot Floer theory. See also~\cite[Proposition~11.57]{LipshitzOzsvathThurston2018} and~\cite[Section~3.2]{DaiHomStoffregenTruong2021}.

\begin{proposition}\label{stru}
    Up to graded chain homotopy equivalence, $C$ is a direct sum of the form
    \[
        \bigoplus_{j=1}^{M}\mathbb F[U]\langle z_j\rangle
        \oplus\bigoplus_{i=1}^{N}
        \bigl(\mathbb F[U]\langle y_i\rangle
        \xrightarrow{U^{\ell_i}}\mathbb F[U]\langle x_i\rangle\bigr),
    \]
    where $\partial z_j=\partial x_i=0$, $\partial y_i=U^{\ell_i}x_i$, and $\ell_i\ge1$. In this model the differential vanishes modulo $U$.
\end{proposition}

\begin{proof}
    Write $R=\mathbb F[U]$ and choose a homogeneous basis of $C$. If
    $\partial=0$, there is nothing to prove. Otherwise, let $\ell$ be the
    smallest exponent of $U$ occurring in the matrix of $\partial$, and write
    $\partial=U^\ell \m{d}$. Then $\m{d}^2=0$, and we may choose a basis element $y$
    such that $x=\m{d}y$ is nonzero modulo $U$.

    Since $\m{d}x=0$, the reduction of $x$ modulo $U$ cannot be a scalar
    multiple of the reduction of $y$. Thus $x$ has a nonzero constant
    coefficient on some basis element other than $y$. By homogeneity, this
    coefficient is a unit, so we may replace that basis element by $x$.
    We now have $\m{d}y=x$ and $\m{d}x=0$.

    For each remaining basis element $w$, let $c_w$ be the coefficient of
    $x$ in $\m{d}w$, and replace $w$ by $w-c_wy$. The differential of each
    resulting generator has zero $x$-coefficient. Its $y$-coefficient is
    also zero, since this is the $x$-coefficient of its image under $\m{d}^2$.
    Hence the remaining generators span a subcomplex, and the pair
    $\partial y=U^\ell x$ splits off.

    Repeating this argument gives the stated decomposition. Pairs with
    $\ell=0$ are contractible and may be removed. All remaining
    differentials vanish modulo $U$.
\end{proof}

In the notation of Proposition~\ref{stru}, we have
\begin{equation}\label{homology}
    \m{H}(C)\cong \bigoplus_{j=1}^M\mathbb{F}[U]\langle z_j\rangle
    \oplus \bigoplus_{i=1}^N \bigl(\mathbb{F}[U]/(U^{\ell_i})\bigr)\langle x_i\rangle.
\end{equation}
Thus $\m{Ord}_U(C)=\max(\{0\}\cup\{\ell_i\}_{i=1}^N)$.

The short exact sequence $0\to C\xrightarrow{U}C\xrightarrow{j}\widehat{C}\to 0$ induces the long exact sequence
\begin{equation}\label{pleq}
    \cdots\longrightarrow\m{H}_{\delta+1}(C)\xrightarrow{U}
    \m{H}_{\delta}(C)\xrightarrow{j_*}\widehat{V}_{\delta}
    \xrightarrow{k_*}\m{H}_{\delta}(C)\xrightarrow{U}
    \m{H}_{\delta-1}(C)\longrightarrow\cdots,
\end{equation}
where $\widehat{V}_\delta=\bigoplus_a \widehat{V}_{\delta,a}$. The map $\m{d}_1=j_*k_*$ preserves $\delta$, and $\m{d}_1^2=0$ since $k_*j_*=0$. This construction is natural under $\mathbb F[U]$-linear chain homotopy equivalences. We may therefore compute $\m{d}_1$ in the model of Proposition~\ref{stru}, keeping the Alexander grading on $\widehat{V}$.

\begin{proposition}\label{prop4.2}
    Suppose that $\m{d}_1\bigl(\widehat{V}_{\delta,a}\bigr)\subset \widehat{V}_{\delta,a-1}\oplus \widehat{V}_{\delta,a+1}$. Then for every $\delta$,
    \begin{equation}\label{keyineq}
        \dim \widehat{V}_{\delta,g}-\dim \widehat{V}_{\delta,g-1}
        \le \dfrac{1}{2}\dim \m{H}\bigl(\widehat{V},\m{d}_1\bigr).
    \end{equation}
\end{proposition}

\begin{proof}
Fix an integer $\delta$. Since $\m{d}_1$ shifts the Alexander grading
by $\pm1$ and $\widehat{V}$ is supported in Alexander gradings $[-g,g]$, we have
\begin{align*}
    \m{d}_{1,g}
    :=\m{d}_1|_{\widehat{V}_{\delta,g}\oplus \widehat{V}_{\delta,-g}}
    &\colon \widehat{V}_{\delta,g}\oplus \widehat{V}_{\delta,-g}
    \longrightarrow \widehat{V}_{\delta,g-1}+\widehat{V}_{\delta,-g+1},\\
    \m{d}_{1,g-1}
    :=\m{d}_1|_{\widehat{V}_{\delta,g-1}+\widehat{V}_{\delta,-g+1}}
    &\colon \widehat{V}_{\delta,g-1}+\widehat{V}_{\delta,-g+1}
    \longrightarrow \widehat{V}_\delta.
\end{align*}
Here the sum $\widehat{V}_{\delta,g-1}+\widehat{V}_{\delta,-g+1}$ is taken inside
$\widehat{V}_\delta$, so it is a single copy of $\widehat{V}_{\delta,0}$ when $g=1$.
Composing $\m{d}_{1,g-1}$ with the projection
$\Pi_g\colon \widehat{V}_\delta\to \widehat{V}_{\delta,g}\oplus \widehat{V}_{\delta,-g}$
gives a map back to the two extreme gradings.

Since $\m{d}_1^2=0$, we have
\[
    (\Pi_g\circ\m{d}_{1,g-1})\circ\m{d}_{1,g}
    =\Pi_g\circ\m{d}_1^2|_{\widehat{V}_{\delta,g}\oplus \widehat{V}_{\delta,-g}}=0.
\]
Hence $\m{Im}\,\m{d}_{1,g}\subseteq
\ker(\Pi_g\circ\m{d}_{1,g-1})$. By  the symmetry
$\dim \widehat{V}_{\delta,a}=\dim \widehat{V}_{\delta,-a}$, we obtain
\begin{align*}
    2\dim \widehat{V}_{\delta,g-1}-\m{rank}(\Pi_g\circ\m{d}_{1,g-1})
    &\ge \dim\ker(\Pi_g\circ\m{d}_{1,g-1})\\
    &\ge \m{rank}\,\m{d}_{1,g}
     =2\dim \widehat{V}_{\delta,g}-\dim\ker\m{d}_{1,g}.
\end{align*}
Thus
\begin{equation}\label{4.3}
    2\bigl(\dim \widehat{V}_{\delta,g}-\dim \widehat{V}_{\delta,g-1}\bigr)
    \le \dim\ker\m{d}_{1,g}
       -\m{rank}(\Pi_g\circ\m{d}_{1,g-1}).
\end{equation}

Only $\widehat{V}_{\delta,g-1}+\widehat{V}_{\delta,-g+1}$ can contribute to the
$\widehat{V}_{\delta,g}\oplus \widehat{V}_{\delta,-g}$ component of a boundary. Hence
\[
    \ker\m{d}_{1,g}\cap\m{Im}\,\m{d}_1\subseteq \bigl(\widehat{V}_{\delta,g}\oplus \widehat{V}_{\delta,-g}\bigr)\cap \m{Im}\,\m{d}_1 
    \subseteq \m{Im}(\Pi_g\circ\m{d}_{1,g-1}).
\]
Since $\ker\m{d}_{1,g}\subseteq\ker\m{d}_1$, it follows that
\begin{equation}\label{4.5}
    \begin{aligned}
        \dim\ker\m{d}_{1,g}-\m{rank}(\Pi_g\circ\m{d}_{1,g-1})
        &\le \dim\frac{\ker\m{d}_{1,g}}
                       {\ker\m{d}_{1,g}\cap\m{Im}\,\m{d}_1}\\
        &\le \dim\frac{\ker\m{d}_1}{\m{Im}\,\m{d}_1}
         =\dim\m{H}\bigl(\widehat{V},\m{d}_1\bigr).
    \end{aligned}
\end{equation}
Combining \eqref{4.3} and \eqref{4.5} gives
\[
    \dim \widehat{V}_{\delta,g}-\dim \widehat{V}_{\delta,g-1}
    \le \frac12\dim\m{H}\bigl(\widehat{V},\m{d}_1\bigr),
\]
as desired.
\end{proof}

\begin{proposition}\label{prop4.3}
    In the notation of Proposition~\ref{stru}, there is a $\delta$-graded isomorphism
    \begin{equation}\label{Hcd1}
        \m{H}\bigl(\widehat{V},\m{d}_1\bigr)\cong\bigoplus_{j=1}^M\mathbb F\langle z_j\rangle
        \oplus\bigoplus_{\ell_i\ge2}\mathbb F\langle x_i,y_i\rangle.
    \end{equation}
\end{proposition}

\begin{proof}
Work in the model of Proposition~\ref{stru}. Since the differential
vanishes modulo $U$, we identify $\widehat{V}$ with $\widehat{C}$ and use the same
symbols for the basis elements of $C$ and their reductions modulo $U$.
For a cycle $x\in C$, the map $j_*$ is given by
$j_*([x])=x\bmod U$. For $x\in\widehat{C}$, the connecting homomorphism is
given by $k_*(x)=[U^{-1}\partial x]$, where on the right we use a lift
of $x$ to $C$. Thus
\begin{align*}
    \m{d}_1 z_j
    &=j_*k_*(z_j)=j_*([U^{-1}\partial z_j])=0,\\
    \m{d}_1 x_i
    &=j_*k_*(x_i)=j_*([U^{-1}\partial x_i])=0,\\
    \m{d}_1 y_i
    &=j_*k_*(y_i)=j_*([U^{\ell_i-1}x_i])
      =\begin{cases}
          x_i & \text{if } \ell_i=1,\\
          0   & \text{if } \ell_i\ge 2.
       \end{cases}
\end{align*}
It follows that
\begin{align}
    \ker\m{d}_1
    &=\bigoplus_{j=1}^M\mathbb{F}\langle z_j\rangle
      \oplus\bigoplus_{i=1}^N\mathbb{F}\langle x_i\rangle
      \oplus\bigoplus_{\ell_i\ge 2}\mathbb{F}\langle y_i\rangle,\\
    \m{Im}\,\m{d}_1
    &=\bigoplus_{\ell_i=1}\mathbb{F}\langle x_i\rangle.\label{imd1}
\end{align}
Taking the quotient gives the $\delta$-graded isomorphism
\[
    \m{H}\bigl(\widehat{V},\m{d}_1\bigr)
    =\frac{\ker\m{d}_1}{\m{Im}\,\m{d}_1}
    \cong\bigoplus_{j=1}^M\mathbb{F}\langle z_j\rangle
    \oplus\bigoplus_{\ell_i\ge 2}\mathbb{F}\langle x_i,y_i\rangle,
\]
as desired.
\end{proof}

\begin{theorem}\label{algorderone}
    Suppose that $\m{d}_1\bigl(\widehat{V}_{\delta,a}\bigr)\subset \widehat{V}_{\delta,a-1}\oplus \widehat{V}_{\delta,a+1}$, that $\m{Ord}_U(C)=1$, and that $\m{H}(C)$ has a single free $\mathbb F[U]$-summand. Then for every integer $\delta$,
    \[
        \dim \widehat{V}_{\delta,g-1}\ge\dim \widehat{V}_{\delta,g}.
    \]
\end{theorem}

\begin{proof}
    By \eqref{homology}, the number of the free summands of $\m{H}(C)$  is $M$, so $M=1$. Since $\m{Ord}_U(C)=1$, every torsion exponent $\ell_i$ equals one. Proposition~\ref{prop4.3} now gives $\dim\m{H}\bigl(\widehat{V},\m{d}_1\bigr)=1$. In particular, \eqref{keyineq} yields
    \[
        2\bigl(\dim \widehat{V}_{\delta,g}-\dim \widehat{V}_{\delta,g-1}\bigr)
        \le \dim\m{H}\bigl(\widehat{V},\m{d}_1\bigr)\le 1.
    \]
    The left-hand side is even, so it is at most zero.
\end{proof}

\begin{corollary}\label{general}
    Let $K\subset S^3$ be a knot of Seifert genus $g$. If $\m{Ord}_U(K)=1$ or $\m{Ord}'_U(K)=1$, then for every integer $\delta$,
    \[
        \dim\hfk_{g-1+\delta}(K,g-1)\ge \dim\hfk_{g+\delta}(K,g).
    \]
\end{corollary}

\begin{proof}
Take $C=\cfk^-(K)$ in the oriented case and $C=\cfk'(K)$ in the
unoriented case. In both cases, $\m{H}(C)$ has a single free $\mathbb F[U]$-summand, and $\widehat{V}_{\delta,a}=\hfk_{a+\delta}(K,a)$. The torsion
order hypothesis excludes the unknot, so $g>0$. By genus detection and
knot Floer symmetry~\cite{OzsvathSzabo2004b,OzsvathSzabo2004},
$\widehat{V}_{\delta,a}=0$ for $|a|>g$ and
$\dim \widehat{V}_{\delta,a}=\dim \widehat{V}_{\delta,-a}$.

It remains to check how $\m{d}_1$ changes the Alexander grading.
On the hat complex, consider the maps
\begin{align*}
    \Phi(\mathbf{x})
    &=\sum_{\mathbf{y}\in\mathbb{T}_{\bm{\alpha}}\cap\mathbb{T}_{\bm{\beta}}}
      \ \sum_{\substack{\phi\in\pi_2(\mathbf{x},\mathbf{y}),\ \mu(\phi)=1\\
                        n_w(\phi)=1,\ n_z(\phi)=0}}
      \#\widehat{\mathcal{M}}(\phi)\,\mathbf{y},\\
    \Psi(\mathbf{x})
    &=\sum_{\mathbf{y}\in\mathbb{T}_{\bm{\alpha}}\cap\mathbb{T}_{\bm{\beta}}}
      \ \sum_{\substack{\phi\in\pi_2(\mathbf{x},\mathbf{y}),\ \mu(\phi)=1\\
                        n_w(\phi)=0,\ n_z(\phi)=1}}
      \#\widehat{\mathcal{M}}(\phi)\,\mathbf{y}.
\end{align*}
These are the coefficients of $U$ and $V$, respectively, in the
differential of $\ccfl(K)$. Comparing the linear terms in $\partial^2=0$
shows that both maps commute with the hat differential. We use the same
symbols for the induced maps on $\hfk(K)$. Their $(M,A)$-degree shifts
are $(1,1)$ and $(-1,-1)$, respectively, so both preserve $\delta=M-A$.

For a cycle $x$ in the hat complex, take its constant-coefficient lift
to $C$. By the definition of the connecting homomorphism,
\[
    \m{d}_1[x]=j_*k_*[x]
    =\bigl[(U^{-1}\partial x)\bmod U\bigr].
\]
Thus setting $V=0$ gives $\m{d}_1^-=\Phi$, while setting $V=U$ gives
$\m{d}_1'=\Phi+\Psi$. In either case,
\[
    \m{d}_1\bigl(\widehat{V}_{\delta,a}\bigr)
    \subseteq \widehat{V}_{\delta,a-1}\oplus \widehat{V}_{\delta,a+1}.
\]
Theorem~\ref{algorderone} now gives
$\dim \widehat{V}_{\delta,g-1}\ge\dim \widehat{V}_{\delta,g}$, as desired.
\end{proof}

Taking $\delta=m-g$ in Corollary \ref{general} recovers Corollary~\ref{tor1} and proves Theorem~\ref{tor}.

\appendix

\section{The relation between the two torsion orders}\label{sec:torsion-orders}

Theorems~\ref{ord} and~\ref{tor} use the oriented and unoriented torsion orders, respectively, to give different bounds on the next-to-top term of knot Floer homology. This naturally leads us to ask how these two torsion orders are related.

\begin{proposition}\label{ord1tuiord'1}
    If a knot $K\subset S^3$ satisfies $\m{Ord}_U(K)=1$, then $\m{Ord}'_U(K)=1$.
\end{proposition}

\begin{proof}
Take $C^-=\m{CFK}^-(K)$ and $C'=\m{CFK}'(K)$. As in the proof of Corollary~\ref{general}, both complexes satisfy the hypotheses of Proposition~\ref{stru}, and their homologies have  a single free $\mathbb F[U]$-summand~\cite{OzsvathSzabo2004,OzsvathStipsiczSzabo2017}. Thus, up to graded chain homotopy equivalence, we have
\begin{equation}\label{algminus}
    \begin{aligned}
        C^- &\simeq \mathbb F[U]\langle z\rangle \oplus \bigoplus_{i=1}^N
        \bigl(\mathbb F[U]\langle x_i\rangle \xrightarrow{U^{n_i}} \mathbb F[U]\langle y_i\rangle\bigr),\\
        C' &\simeq \mathbb F[U]\langle z'\rangle \oplus \bigoplus_{j=1}^{N'}
        \bigl(\mathbb F[U]\langle x'_j\rangle \xrightarrow{U^{\ell_j}} \mathbb F[U]\langle y'_j\rangle\bigr).
    \end{aligned}
\end{equation}
Since $C^-/U\cong C'/U\cong\widehat{C}$, the two reductions have isomorphic homology. The differentials in \eqref{algminus} vanish modulo $U$, so $\dim\m{H}(C^-/U)=1+2N$ and $\dim\m{H}(C'/U)=1+2N'$. Thus $N=N'$.

As shown in the proof of Corollary~\ref{general}, $\m{d}_1^-=\Phi$ and $\m{d}_1'=\Phi+\Psi$. Since $\Phi$ is homogeneous in the Alexander grading, choose homogeneous elements $w_1,\ldots,w_r$ such that $\Phi w_1,\ldots,\Phi w_r$ form a basis of $\m{Im}\,\Phi$. By \eqref{imd1} and the hypothesis $\m{Ord}_U(K)=1$, we have $r=\m{rank}\,\Phi=\#\{i:n_i=1\}=N$. We claim that $(\Phi+\Psi)w_1,\ldots,(\Phi+\Psi)w_r$ are also linearly independent. Suppose that
\begin{equation}\label{linear}
    c_1(\Phi+\Psi)w_1+\cdots+c_r(\Phi+\Psi)w_r=0.
\end{equation}
Let $a=\displaystyle\max_{1\le k\le r} A(w_k)$ and reorder the $w_k$ so that $w_1,\ldots,w_{k_1}$ are precisely those of Alexander grading $a$. Since $\Phi$ and $\Psi$ shift the Alexander grading by $+1$ and $-1$, respectively, the component of \eqref{linear} in grading $a+1$ is $c_1\Phi w_1+\cdots+c_{k_1}\Phi w_{k_1}=0$. By linear independence, $c_1=\cdots=c_{k_1}=0$. Repeating the argument in decreasing order of Alexander grading gives $c_k=0$ for every $k$. Thus $\m{rank}\,(\Phi+\Psi)\ge\m{rank}\,\Phi=N=N'$.

On the other hand, \eqref{imd1} and \eqref{algminus} give $\m{rank}\,(\Phi+\Psi)=\#\{j:\ell_j=1\}\le N'$. Hence $\#\{j:\ell_j=1\}=N'$, so every $\ell_j$ equals one. Hence we conclude that $\m{Ord}'_U(K)=1$.
\end{proof}

In general, Himeno and Teragaito asked which pairs of torsion orders can occur.

\begin{question}[{\cite[Question~1.4]{HimenoTeragaito2024}}]
    Given positive integers $m$ and $n$, is there a knot $K\subset S^3$ such that $\bigl(\m{Ord}_U(K),\m{Ord}'_U(K)\bigr)=(m,n)$?
\end{question}

Proposition~\ref{ord1tuiord'1} rules out the pairs $(1,n)$ for $n\ge 2$. For each fixed $n$, the following proposition shows that there exists a lower bound for the possible values of \(m\).

\begin{proposition}
    For any positive integer $n$, there exists an integer $M(n)$ such that for any positive integer $m$, there exists a knot $K$ such that $\bigl(\m{Ord}_U(K),\m{Ord}'_U(K)\bigr)=(m,n)$ if and only if $m\ge M(n)$. Moreover, we have $M(1)=1$, $M(2)=2$ and $2\le M(n)\le 2n-1$ for $n\ge 3$.
\end{proposition}

\begin{proof}
We first construct a knot $J_m$ realizing $(m,1)$ for each $m\ge 1$. Recall that a positive torus knot $T_{p,q}$ with $1<p<q$ has bridge index $p$ and oriented torsion order $p-1$~\cite[Corollary~5.3]{JuhaszMillerZemke2020}. Thus the trefoil $J_1=T_{2,3}$ realizes $(1,1)$ by Proposition~\ref{ord1tuiord'1}. For even \(m\ge 2\), take \(J_m=T_{m+1,m+3}\), which has oriented torsion order \(m\). A single rational tangle replacement, specifically a pinch move~\cite[Figure~1, Lemma~2.1 and Proposition~5.8]{Himeno2025band}, takes \(J_m\) to the unknot \(T_{1,1}\); hence, by the proof of Proposition~\ref{3}, \(J_m\) has unoriented torsion order one.

For odd $m=2k-1\ge 3$, let $J_m$ be the $(2,1)$-cable of $T_{k,2k+1}$. The companion is an L-space knot distinct from the trefoil, and its bridge index is $k=\m{Ord}_U(T_{k,2k+1})+1$. The equality case of~\cite[Proposition~1.4]{HomLidmanPark2026} therefore gives
\[
    \m{Ord}_U(J_m)=2\bigl(\m{Ord}_U(T_{k,2k+1})+1\bigr)-1=2k-1=m.
\]
The rational tangle replacement in Figure~\ref{fig:solid-torus-patterns} unknots $J_m$. Since $J_m$ is nontrivial, its rational unknotting number is one, so $\m{Ord}'_U(J_m)=1$ by the proof of Proposition~\ref{3}.

Now fix $n\ge 1$. The torus knot $T_{2n,2n+1}$ realizes $(2n-1,n)$~\cite[Theorem~1.3]{Eftekhary2025rational}. Let $M(n)$ be the least positive integer for which $(M(n),n)$ is realized, and choose a knot $K$ realizing this pair. For any $m\ge M(n)$, both torsion orders of $K\#J_m$ are given by the connected sum formula~\cite[Proposition~1.6] {Eftekhary2025rational} $$ \m{Ord}_U(K\#K')=\max\{ \m{Ord}_U(K),\m{Ord}_U(K')\},\quad \m{Ord}_U'(K\#K')=\max\{ \m{Ord}_U'(K),\m{Ord}_U'(K')\},$$hence we have 
$$  \bigl(\m{Ord}_U(K\#J_m),\m{Ord}'_U(K\#J_m)\bigr)
    =\bigl(\max\{M(n),m\},\max\{n,1\}\bigr)=(m,n).$$
This proves the claimed equivalence.

The trefoil gives $M(1)=1$. For $n\ge 2$, Proposition~\ref{ord1tuiord'1} gives $M(n)\ge 2$, while $T_{2n,2n+1}$ gives $M(n)\le 2n-1$. Finally, $12n_{404}$ realizes $(2,2)$~\cite[Example~5.2]{Eftekhary2025rational}, so $M(2)=2$.
\end{proof}

\begin{remark}
The preceding proposition answers Himeno and Teragaito's question for \(n=1\) and \(n=2\). For \(n\ge 3\), it remains to determine \(M(n)\), which we leave for future work.
\end{remark}

\end{document}